\documentclass{amsart}
\usepackage[utf8]{inputenc}
\usepackage{amsfonts}
\usepackage{amssymb,amsmath,amscd,euscript,verbatim,array}
\usepackage{geometry}
\usepackage{indentfirst}
\usepackage{graphicx}
\usepackage{ifpdf}
\usepackage{enumitem}
\usepackage{mathrsfs}

\usepackage{bm}
\usepackage[hidelinks]{hyperref}
\makeatletter
\renewcommand{\l@subsection}{\@tocline{2}{0pt}{2.5pc}{5pc}{}}
\makeatother
\newcommand{\T}{\mathbb T}
\newcommand{\R}{\mathbb R}

\newcommand{\Z}{\mathbb Z}

\newtheorem{Theorem}{Theorem}
\newtheorem{Lemma}{Lemma}[section]
\newtheorem{Proposition}[Lemma]{Proposition}
\theoremstyle{definition}
\newtheorem{Definition}[Lemma]{Definition}
\newtheorem{Remark}[Lemma]{Remark}
\begin{document}
\title[Persistence of invariant graphs for twist maps]{Persistence of invariant graphs for twist maps under analytic perturbations}
\author{Qi Li}
\address{School of Mathematics and Statistics, Beijing Institute of Technology, Beijing 100081, China}
\email{qilicindy@bit.edu.cn}
\author{Yi Liu}
\address{School of Mathematics and Statistics, Beijing Institute of Technology, Beijing 100081, China}
\email{yiliu111@foxmail.com}
\author{Lin Wang}
\address{School of Mathematics and Statistics, Beijing Institute of Technology, Beijing 100081, China}
\email{lwang@bit.edu.cn}
\begin{abstract}
We construct high-frequency perturbations of an integrable exact
symplectic twist map which preserve an analytic invariant graph with
a prescribed rotation number. Their oscillation frequencies may be
any positive integers tending to infinity. For every fixed Gevrey
parameter $L>0$, the perturbations
tend to zero in a Gevrey-$\alpha$ norm and diverge in a
Gevrey-$(\alpha-\varepsilon)$ norm. The frequency class strictly
contains the $\alpha$-Brjuno--R\"ussmann class, and  their difference has cardinality \(2^{\aleph_0}\).
\end{abstract}
\maketitle
\tableofcontents
\clearpage

\section{Introduction and main result}

We study analytic invariant graphs of exact symplectic twist maps
whose dynamics is analytically conjugate to a prescribed irrational
rotation. We construct increasingly oscillatory perturbations that
preserve such graphs and have different asymptotic sizes in two
Gevrey norms.

\subsection{Background and setting}

Moser's
invariant-curve theorem~\cite{Mos62} established their persistence
for sufficiently small smooth perturbations under a Diophantine
condition. The circle-diffeomorphism results of Herman~\cite{Her79}
and Yoccoz~\cite{Yoc95} describe how arithmetic affects
the regularity of conjugacies to rotations. In the construction below,
the perturbing functions belong to a Gevrey-$\gamma$ class with
$0\le\gamma<1$. The range $0<\gamma<1$ consists of ultra-analytic
functions; the endpoint $\gamma=0$ will mean trigonometric polynomials.
We use this regularity to estimate analytic norms, while measuring
the size of the perturbations in Gevrey-$\alpha$ norms with $\alpha>1$.
In the Hamiltonian setting, Popov~\cite{Pop04} proved a Gevrey
KAM theorem with Whitney--Gevrey dependence of the invariant tori
on Diophantine frequencies.

Write $\T:=\R/\Z$. Let $F=(X,Y):\R^2\to\R^2$ be a
diffeomorphism satisfying
\[
 F(x+1,y)=F(x,y)+(1,0).
\]
It induces a map $f$ of $\T\times\R$. We call $f$ an exact
symplectic twist map if it is isotopic to the identity, if
$\partial_yX>0$ and $(x,y)\mapsto(x,X(x,y))$ is a diffeomorphism,
and if there is a generating function $S(x,X)$ satisfying
\[
 S(x+1,X+1)=S(x,X),\qquad
 Y\,dX-y\,dx=dS(x,X).
\]
The differential identity is evaluated on the graph of $F$;
equivalently,
\[
 y=-\partial_1S(x,X),\qquad Y=\partial_2S(x,X).
\]
Taking exterior derivatives gives $dX\wedge dY=dx\wedge dy$,
so exact symplecticity implies area preservation. We use the same
notation for a cylinder map and its lift when no ambiguity arises.
\begin{Definition}[\( C^0 \)-invariant graph]
    \label{def:invariant_graph}
    \( \Gamma \subset \mathbb{T} \times \mathbb{R} \) is called a \( C^0 \)-invariant graph of \( f \) if
    \begin{enumerate}
        \item[(i)] \( \Gamma = \{ (x, \Psi(x)) : x \in \mathbb{T} \}, \) where \( \Psi : \mathbb{T} \to \mathbb{R} \) is a continuous function;
        \item[(ii)] \( \Gamma \) is invariant under the action of \( f \).
    \end{enumerate}
\end{Definition}

\begin{Remark}\label{rem:properties}
\
    \begin{enumerate}
        \item[(i)] Under the area-preserving and twist conditions, Birkhoff proved that if \( f \) is of class \( C^1 \), then \( \Psi \) is a Lipschitz function on \( \mathbb{T} \) (see \cite[Chapter~I, Proposition~2.2]{H1}). Moreover, for area-preserving twist maps, Birkhoff also proved that if an essential invariant circle (i.e., a noncontractible invariant Jordan curve) exists, then it must be the graph of a function. These two results are collectively referred to as the Birkhoff graph theorem.
        \item[(ii)] The invariance condition can equivalently be stated for the lift: if \( \tilde{\Psi} : \mathbb{R} \to \mathbb{R} \) is the 1-periodic lift of \( \Psi \), then the graph \( \tilde{\Gamma} := \{ (x, \tilde{\Psi}(x)) \mid x \in \mathbb{R} \} \) is invariant under \( F \).
    \end{enumerate}
\end{Remark}

The restriction to an invariant graph is an orientation-preserving
circle homeomorphism. We require the graph and its conjugacy to
a rotation to be analytic, as specified below.

\begin{Definition}\label{def:property}
An invariant graph has Property $A_\omega$ if its defining function
$\Psi$ is real analytic and the dynamics on the graph is
real-analytically conjugate to the rigid rotation
$R_\omega(\theta):=\theta+\omega$ on $\T$.
\end{Definition}

Consider the map
\begin{equation}\label{basmap}
 f(x,y)=(x+\omega+y+\eta\phi(x),\,y+\eta\phi(x)),
\end{equation}
where $\omega\in(0,1)\setminus\mathbb Q$, $\eta\in\R$, and
$\phi$ is real analytic, one-periodic, and has zero average. This is a
perturbation of the integrable map
\[
 T_\omega(x,y):=(x+\omega+y,y),
\]
whose graph $y=0$ has Property $A_\omega$.

The results above motivate the following quantitative question.
\vspace{1ex}
\begin{itemize}
    \item \textbf{Question:} \emph{If the perturbation \(\phi\) has higher regularity than that required by the topology used to measure its size, can the invariant graph with Property \(A_\omega\) be preserved for a larger class of frequencies \(\omega\)?}
\vspace{1ex}

\end{itemize}
In \cite[Corollary 4]{LW25}, the authors made a preliminary exploration of this question by constructing explicit examples. In this paper, we are devoted to studying the mechanism behind the question and providing results that are as general as possible. More precisely, we seek perturbations that are small in a Gevrey norm and preserve analytic invariant graphs for rotation numbers outside the corresponding Gevrey arithmetic class. Smallness in a
Gevrey-$\alpha$ norm, with $\alpha>1$, does not imply smallness
in an analytic norm on a fixed strip. We instead use ultra-analytic
Fourier decay to estimate the analytic norm on a strip whose width
depends on the oscillation frequency. The rotation number
$\omega$ remains fixed.

\subsubsection{Regularity and arithmetic preliminaries}
For $\gamma>0$ and $L>0$, define
\begin{equation}\label{gevr}
 \|\phi\|_{\gamma,L}:=
 \sum_{k=0}^{\infty}\frac{L^{k\gamma}}{(k!)^\gamma}
 \|\partial^k\phi\|_{C^0(\T)}.
\end{equation}
Write
\[
 G^{\gamma,L}(\T):=\{\phi\in C^\infty(\T):
                         \|\phi\|_{\gamma,L}<\infty\},\qquad
 G^\gamma(\T):=\bigcup_{L>0}G^{\gamma,L}(\T).
\]
We use this derivative-norm convention as in \cite{MS03};
$\gamma=1$ is the analytic case, and $0<\gamma<1$ gives
ultra-analytic classes. The index zero is a separate convention:
\[
 G^0(\T):=\left\{\phi(x)=\sum_{|j|\leq N}
       \widehat\phi(j)e^{2\pi ijx}:N\in\Z_{\ge0}\right\}.
\]
We do not substitute $\gamma=0$ into \eqref{gevr}.
For a bounded vector-valued perturbation $U$ on $\T\times\R$, set
\[
 \|U\|_{\gamma,L}:=
 \sum_{k=0}^{\infty}\frac{L^{k\gamma}}{(k!)^\gamma}
 \max_{i+j=k}\|\partial_x^i\partial_y^j U\|_\infty,
\]
where the vector norm is the maximum of the component norms.
In particular, $\|(\phi(x),\phi(x))\|_{\gamma,L}
=\|\phi\|_{\gamma,L}$. The topology in the main theorem is
Gevrey-$\alpha$, with $\alpha>1$.
Throughout the paper, $\langle g\rangle:=\int_0^1g(\theta)\,d\theta$
denotes the average of a periodic function.

Write $\omega=[0;a_1,a_2,\ldots]$ and let $p_n/q_n$ be its
convergents, with $p_{-1}=1$, $p_0=0$, $q_{-1}=0$, and $q_0=1$.
Thus $q_{n+1}=a_{n+1}q_n+q_{n-1}$ and $q_{n+2}\ge2q_n$.
Write $\|t\|:=\min_{\ell\in\Z}|t-\ell|$ for the distance to the
nearest integer; the notation $\|t\|_{\R/\Z}$ has the same meaning.
For $\alpha>1$, define the Brjuno--R\"ussmann class
\[
 \mathcal{BR}_\alpha:=
 \left\{\omega:\sum_{n\ge0}
                  \frac{\log q_{n+1}}{q_n^{1/\alpha}}<\infty\right\}.
\]
The same formula at index one defines the Brjuno class
$\mathcal{BR}_1$. For integers $k\ge1$, introduce
\[
 \Delta_k^\alpha:=\left\{\omega:
 \sum_{n\ge m}\frac{\log q_{n+1}}{q_n}
       \le k q_m^{1/\alpha-1}\text{ for every }m\ge0\right\},
 \qquad \Delta^\alpha:=\bigcup_{k\ge1}\Delta_k^\alpha.
\]
All frequency sets consist of irrationals in $(0,1)$. In
Section~\ref{sec:arithm} we prove that
\[
 \omega\in\Delta^\alpha
 \quad\Longleftrightarrow\quad
 \log q_{n+1}=O(q_n^{1/\alpha}),
\]
and that
\[
 \mathcal{BR}_\alpha\subsetneq\Delta^\alpha\subseteq\mathcal{BR}_1.
\]
Both \(\mathcal{BR}_\alpha\) and \(\Delta^\alpha\) have full Lebesgue measure and Hausdorff dimension one. Their difference has cardinality \(2^{\aleph_0}\), Lebesgue measure zero, and Hausdorff dimension zero.
Thus the frequency class used below is strictly larger than
$\mathcal{BR}_\alpha$, while still consisting of Brjuno numbers.

\subsection{Main result}

Analytic invariant curves with a prescribed Brjuno rotation number
persist for sufficiently small analytic perturbations of an
integrable exact symplectic twist map. Gentile~\cite{Ge1} gives a
proof by multiscale analysis of the Lindstedt series. In the
Hamiltonian setting, R\"ussmann~\cite{Rus01} developed analytic
invariant-torus theorems under arithmetic conditions weaker than
the Diophantine condition. Bounemoura and F\'ejoz~\cite{BF}
prove a KAM theorem under the corresponding
$\alpha$-Brjuno--R\"ussmann condition for Gevrey perturbations,
including an exact symplectic map version in their Theorem~F.
The distinction between persistence of an invariant graph and analytic
linearizability of its dynamics is relevant to questions of arithmetic
optimality; see Bounemoura~\cite{Bou}. We use ultra-analytic
Fourier bounds to choose explicit amplitudes for a family of
rescaled perturbations. The resulting maps preserve analytic
invariant graphs, while their perturbations tend to zero in one
Gevrey norm and diverge in a stronger one.

Fix $\omega\in\Delta^\alpha$ and let $(N_m)_{m\ge1}$ be any
sequence of positive integers with $N_m\to\infty$. The numbers
$N_m$ specify the oscillation frequencies of the perturbations;
the notation $q_n$ continues to denote the convergent denominators
of $\omega$. No relation between the two sequences is assumed.
For $\varphi\in G^\gamma(\T)$, put
\[
 \phi_m(x):=\varphi(N_mx).
\]
Consider the sequence of twist maps
\begin{equation}\label{yp}
 f_m(x,y):=(x+\omega+y+\delta_m\phi_m(x),\,
                            y+\delta_m\phi_m(x)).
\end{equation}
The derivative of the first component with respect to $y$ and the
Jacobian determinant are both one. Since $N_m$ is a positive
integer, $\phi_m$ is one-periodic and
\[
 \int_0^1\phi_m(x)\,dx
 =\frac1{N_m}\int_0^{N_m}\varphi(t)\,dt
 =\int_0^1\varphi(t)\,dt.
\]
If $\varphi$ has zero average, define
$V_\varphi(x):=\int_0^x\varphi(t)\,dt$. Then $V_\varphi$ is
one-periodic, and
\[
 W_m(x):=\frac1{N_m}V_\varphi(N_mx)
\]
is a one-periodic primitive of $\phi_m$. The generating function
\[
 S_m(x,x'):=\tfrac12(x'-x-\omega)^2+\delta_m W_m(x),
 \qquad y=-\partial_xS_m,\quad y'=\partial_{x'}S_m,
\]
verifies exact symplecticity. The inverse lift is given by
\[
 x=x'-\omega-y',\qquad y=y'-\delta_m\phi_m(x'-\omega-y').
\]
 Moreover,
\[
 \|f_m-T_\omega\|_{s,L}
       =\delta_m\|\varphi(N_m\cdot)\|_{s,L}
 \qquad(\delta_m>0).
\]

\begin{Theorem}\label{Th1}
Fix $\alpha>1$, $\omega\in\Delta^\alpha$, and
$0<\varepsilon<\min\{1,\alpha-1\}$. Let $(N_m)_{m\ge1}$ be any
sequence of positive integers tending to infinity. Let
$\varphi\in G^\gamma(\T)$ be real, nonzero, and of zero average, where
$0\le\gamma<\varepsilon/\alpha$. Choose any
\[
 \frac1{\alpha-\gamma}<\tau_*<\frac1{\alpha-\varepsilon},
 \qquad \delta_m:=e^{-N_m^{\tau_*}}.
\]
There is an integer $N_*\ge1$, depending only on
$\alpha,\omega,\gamma,\varphi$ and $\tau_*$, such that, whenever
$N_m\ge N_*$, the map~\eqref{yp} has an invariant graph with
Property $A_\omega$. In particular, this holds for all sufficiently
large $m$. Its conjugacy can be written
$\theta\mapsto(\theta+u_m(\theta),v_m(\theta))$, where
$\langle u_m\rangle=0$ and $\|u_m'\|_{C^0}\to0$.
Moreover, for every fixed $L>0$,
\[
 \|f_m-T_\omega\|_{\alpha,L}\longrightarrow0,
 \qquad
 \|f_m-T_\omega\|_{\alpha-\varepsilon,L}\longrightarrow+\infty.
\]
The threshold $N_*$ is independent of $L$ and of the sequence
$(N_m)$. The same sequence $\delta_m$ works for all fixed $L$.
\end{Theorem}

\begin{Remark}
\label{rem:integer-frequencies}
The oscillation frequencies $(N_m)$ may be chosen independently of the convergent denominators of $\omega$, with no prescribed growth rate; the choice $N_m=q_m$ is included.
In particular, taking $N_m=m$ gives persistence at every sufficiently
large integer oscillation frequency. For a given sequence, the corresponding index $m_0$
such that $N_m\ge N_*$ for all $m\ge m_0$ may depend on that
sequence. The proof of the uniform threshold is given in
Subsection~\ref{sec:newton-main-application}.
\end{Remark}

\begin{Remark}
\label{rem:comparison}
Bounemoura and F\'ejoz~\cite[Theorem~F]{BF} establish persistence
of Gevrey invariant tori, with Gevrey conjugacies to the prescribed
rigid rotation, for sufficiently small Gevrey perturbations of
nondegenerate integrable exact symplectic maps. Their arithmetic
assumption is the $\alpha$-Brjuno--R\"ussmann condition; for a
one-dimensional rotation it corresponds to the class
$\mathcal{BR}_\alpha$ defined above. The larger class
$\Delta^\alpha$ in Theorem~\ref{Th1} is used only for the
prescribed family $\delta_m\varphi(N_m\cdot)$, whose stronger
regularity $G^\gamma$, $\gamma<1$, controls its analytic size.
The theorem combines explicit amplitudes, invariant
graphs with Property $A_\omega$, and convergence to zero in one
Gevrey norm with divergence in a stronger one, for every fixed
$L>0$ and a single amplitude sequence. 
\end{Remark}

The proof uses a scalar quasi-Newton iteration.  At an approximate solution, the
identity also contains a term involving the derivative of the
residual. This term is retained in the quadratic error estimate.
The construction is related to the factorization of Su and de la
Llave \cite[Section~4.1]{SdL12}. Here factorization
means expressing the principal part of the linearized equation,
after a change of unknown, as the composition of two first-order
difference operators and multiplication by known functions; see
\eqref{eq:newton-factorization}. The arithmetic assumption bounds
the inverse of a first-order difference operator on zero-average
functions, with analytic width loss $t$, by
$\exp(Ct^{-1/(\alpha-1)})$. The correction is obtained by solving
two first-order cohomological equations. An integral identity gives
the first compatibility condition, and a scalar normalization enforces
the second. Related factorizations are used in
parameterization methods; see \cite{dlLGJV05}. We prove the
factorization, compatibility identities, and convergence estimates.
Taking the analytic width proportional to
$N_m^{-(\alpha-1)/(\alpha-\gamma)}$ bounds both the
small-divisor loss and the analytic growth of
$\varphi(N_m\cdot)$ by $\exp(CN_m^{1/(\alpha-\gamma)})$.

\subsection*{Organization of the paper}

Section~\ref{sec:arithm} recalls the continued-fraction estimates
used throughout the proof and shows that the tail condition defining
$\Delta^\alpha$ is equivalent to
$\log q_{n+1}=O(q_n^{1/\alpha})$. It then establishes the
inclusions between the arithmetic classes, their measure and category
properties, and the cardinality and Hausdorff dimension of
$\Delta^\alpha\setminus\mathcal{BR}_\alpha$.

Section~\ref{sec:analytic-newton} reduces the invariant-graph problem
to a scalar difference equation and proves the analytic existence
estimate in Proposition~\ref{prop:analytic-newton}. The argument
first bounds the inverse of each first-order difference operator
on zero-average analytic functions. It then derives the exact
compatibility identity and the factorization of the linearized
equation, constructs normalized corrections, and proves a quadratic
error estimate. Explicit strip widths and error bounds prove
convergence and control the derivative of the solution.

Section~\ref{sec:proofth1} applies this estimate to
$\delta_m\varphi(N_m\cdot)$. Subsections~\ref{sec:newton-arithmetic}
and~\ref{sec:newton-rescaled-force} translate the arithmetic
assumption into a small-divisor bound and estimate the forcing on
complex strips. Subsection~\ref{sec:newton-main-application} chooses
the strip widths and amplitudes, gives a threshold uniform over all
integer oscillation frequencies, constructs an analytic inverse of
the angle conjugacy, and verifies the invariant-graph property.
Finally, Subsection~\ref{sec:newton-gevrey-asymptotics} proves the
upper estimate in the Gevrey-$\alpha$ norm and a lower estimate
from one nonzero Fourier coefficient in the
Gevrey-$(\alpha-\varepsilon)$ norm, completing the two limits in
Theorem~\ref{Th1}.

\vspace{1em}

\section{Arithmetic of the frequency class}\label{sec:arithm}

We retain the convergents $p_n/q_n$ of $\omega$ and the sets
$\Delta_k^\alpha$ and $\Delta^\alpha$ defined in the introduction.

\begin{Proposition}[{\cite[Chapter~V, \S7]{Her79}}]\label{prop:continued-fraction-estimates}
Let $\omega=[0;a_1,a_2,\ldots]$ and let $p_n/q_n$ have the
initial conditions stated above. Write $\rho_n:=q_n\omega-p_n$.
For every $n\ge0$,
\[
 \frac1{q_{n+1}+q_n}<|\rho_n|<\frac1{q_{n+1}},
 \qquad \rho_n\rho_{n+1}<0.
\]
Moreover, if $q_n\le k<q_{n+1}$ and $k\ge1$, then
\[
 |k\omega-\ell|\ge |\rho_n|\qquad(\ell\in\mathbb Z).
\]
Consequently
\[
 \left|\omega-\frac{p_n}{q_n}\right|
 <\frac1{q_nq_{n+1}}.
\]
\end{Proposition}

We also use the standard description of continued-fraction
cylinders; see \cite[Chapter~V]{Her79}. The irrationals with fixed
initial partial quotients $a_1,\ldots,a_m$ are exactly the irrational
points of the interval with endpoints $p_m/q_m$ and
$(p_m+p_{m-1})/(q_m+q_{m-1})$. Its length is
$1/(q_m(q_m+q_{m-1}))$, which tends to zero as $m\to\infty$.
These cylinders consequently form a basis for the irrational
points of $(0,1)$.

For each fixed $k$, the set $\Delta_k^\alpha$ is nowhere dense in
$(0,1)$. Indeed, let $I\subset(0,1)$ be a nonempty open interval.
Choose a finite continued-fraction cylinder contained in $I$ and
extend its word, if necessary, to end at index $m$. Its denominator
$q_m$ is now fixed. By choosing $a_{m+1}$ sufficiently large we obtain
\[
 \log q_{m+1}>kq_m^{1/\alpha}.
\]
Every irrational in this smaller cylinder violates the defining
inequality for $\Delta_k^\alpha$ at $m$, already at its first
summand. A cylinder consists of all irrationals in an interval
whose endpoints are rational; that interval is
disjoint from $\Delta_k^\alpha$. Every nonempty interval $I$ thus
contains a nonempty subinterval disjoint from $\Delta_k^\alpha$.
Consequently, the closure of $\Delta_k^\alpha$ has empty interior.
Taking the countable union over $k$ proves that $\Delta^\alpha$ is meager.

We next compare the defining tail bound with a pointwise bound on the continued-fraction denominators.

	\begin{Lemma}\label{fs}
For every $d>0$ and every $m\geq 0$,
\begin{equation}
\sum_{n=m}^{\infty} q_n^{-d}
\leq
\frac{2}{1-2^{-d}}\,q_m^{-d}.
\end{equation}
\end{Lemma}
\begin{proof}
The recurrence gives
\[
q_{n+2}\geq 2q_n,
\]
and hence
\[
q_{m+j}\geq 2^{\lfloor j/2\rfloor}q_m.
\]
Summing the resulting geometric series proves the assertion,
including $m=0$ when the first denominator is repeated.
\end{proof}
\begin{Proposition}\label{bca}
The following statements are equivalent:
\begin{enumerate}
    \item $\omega\in\Delta^\alpha$;
    \item there is a finite constant $H$ such that
    \begin{equation}
    \log q_{n+1}\leq Hq_n^{1/\alpha},
    \qquad n\geq 0.\label{he}
    \end{equation}
\end{enumerate}
Moreover,
\begin{equation}\label{arithmetic-inclusions}
\mathcal{BR}_\alpha
\subseteq
\Delta^\alpha
\subseteq
\mathcal{BR}_1.
\end{equation}
\end{Proposition}
\begin{proof}
If the defining tail bound holds with constant $K_\omega$,
retaining its first term and multiplying by $q_m$ gives (\ref{he})
with $H=K_\omega$.
Conversely, set
\[
d=1-\frac{1}{\alpha}>0.
\]
Lemma~\ref{fs} gives
\[
\sum_{n=m}^{\infty}
\frac{\log q_{n+1}}{q_n}
\leq
H\sum_{n=m}^{\infty}q_n^{-d}
\leq
\frac{2H}{1-2^{-d}}\,
q_m^{1/\alpha-1}.
\]
Enlarging the constant to an integer at least one gives membership in
$\Delta^\alpha$.
If $\omega\in\mathcal{BR}_\alpha$, monotonicity of $q_n$ gives
\[
\begin{aligned}
\sum_{n=m}^{\infty}
\frac{\log q_{n+1}}{q_n}
&=
\sum_{n=m}^{\infty}
\frac{\log q_{n+1}}{q_n^{1/\alpha}}
q_n^{1/\alpha-1}
\\
&\leq
q_m^{1/\alpha-1}
\sum_{n=0}^{\infty}
\frac{\log q_{n+1}}{q_n^{1/\alpha}}.
\end{aligned}
\]
This proves the first inclusion.
The second follows by taking $m=0$ in the defining bound for
$\Delta^\alpha$.
\end{proof}
\begin{Proposition}\label{arithmetic-full-measure}
Both $\mathcal{BR}_\alpha$ and $\Delta^\alpha$ have full Lebesgue
measure in $(0,1)$ and Hausdorff dimension one.
\end{Proposition}
\begin{proof}
For each integer $q\geq 2$, let
\[
E_q
=
\bigcup_{p=0}^{q}
\left(
\frac{p}{q}-\frac{1}{q^3},
\frac{p}{q}+\frac{1}{q^3}
\right)
\cap(0,1).
\]
Its Lebesgue measure is at most
\[
2(q+1)q^{-3},
\]
and the sum of these bounds is finite.
Since
\[
 \left|\bigcup_{q\ge Q}E_q\right|
 \le \sum_{q\ge Q}2(q+1)q^{-3}\longrightarrow0,
\]
the set $\bigcap_{Q\ge2}\bigcup_{q\ge Q}E_q$ has measure zero.
Thus, for almost every irrational $\omega$ there is $Q(\omega)$ such that
\[
\left|\omega-\frac{p}{q}\right|
\geq q^{-3}
\]
whenever $q\geq Q(\omega)$ and $p\in\mathbb{Z}$.
For its convergents, Proposition~\ref{prop:continued-fraction-estimates} gives
\[
\left|
\omega-\frac{p_n}{q_n}
\right|
<
\frac{1}{q_nq_{n+1}},
\]
which, together with the preceding lower bound, implies
\[
q_{n+1}<q_n^2
\]
for all sufficiently large $n$.
Consequently,
\[
\sum_{n=0}^{\infty}
\frac{\log q_{n+1}}{q_n^{1/\alpha}}
<\infty.
\]
Indeed, after omitting finitely many terms, its summands are bounded by
\[
\frac{2\log q_n}{q_n^{1/\alpha}},
\]
which is at most a constant times
\[
q_n^{-1/(2\alpha)};
\]
Lemma~\ref{fs} applies.
Thus $\mathcal{BR}_\alpha$ has full measure, and so does
$\Delta^\alpha$ by \eqref{arithmetic-inclusions}.
A subset of an interval with full Lebesgue measure has Hausdorff
dimension one.
\end{proof}
\begin{Proposition}
The set
\[
\Delta^\alpha\setminus\mathcal{BR}_\alpha
\]
has cardinality $2^{\aleph_0}$, Lebesgue measure zero, and Hausdorff
dimension zero.
\end{Proposition}
\begin{proof}
Choose $Q_*>1$ so large that
\[
e^{q^{1/\alpha}}\geq 2q
\]
for every $q\geq Q_*$, and fix an initial continued-fraction segment
whose final denominator $q_N$ is at least $Q_*$.
Given any binary sequence $(b_n)_{n\geq N}$, define the subsequent
partial quotients by
\begin{equation}\label{spq}
a_{n+1}
=
\left\lfloor
\frac{e^{q_n^{1/\alpha}}}{q_n}
\right\rfloor
+b_n,
\qquad
b_n\in\{0,1\},
\qquad
n\geq N.
\end{equation}
These partial quotients are positive.
The denominator recurrence gives
\[
e^{q_n^{1/\alpha}}-q_n
<
q_{n+1}
<
e^{q_n^{1/\alpha}}+2q_n.
\]
Therefore,
\[
\frac12 e^{q_n^{1/\alpha}}
<
q_{n+1}
<
2e^{q_n^{1/\alpha}},
\qquad
\left|
\log q_{n+1}-q_n^{1/\alpha}
\right|
\leq\log 2.
\]
Since $q_n\to\infty$, it follows that
\begin{equation}
\frac{\log q_{n+1}}{q_n^{1/\alpha}}\label{fl}
\longrightarrow 1.
\end{equation}
The ratios in \eqref{fl} are bounded, so Proposition~\ref{bca}
implies that the resulting irrational number belongs to
$\Delta^\alpha$.
Their limit is nonzero, so their sum diverges and the number does not
belong to $\mathcal{BR}_\alpha$.
If two binary sequences first differ at index $n$, the preceding
partial quotients and $q_n$ agree, whereas (\ref{spq}) gives different
values of $a_{n+1}$.
Uniqueness of the infinite continued-fraction expansion of an
irrational number makes this construction injective.
It produces $2^{\aleph_0}$ distinct numbers; the reverse cardinality
bound follows from containment in $\mathbb{R}$.
The difference set has measure zero because $\mathcal{BR}_\alpha$
has full measure by Proposition~\ref{arithmetic-full-measure}.
For the dimension assertion, for each integer $\tau\geq3$ set
\[
\mathcal{L}_\tau
=
\left\{
x\in(0,1):
\left|x-\frac{p}{q}\right|
<
q^{-\tau}
\text{ for infinitely many }
(p,q)\in\mathbb{Z}\times\mathbb{N}
\right\}.
\]
An irrational $x\notin\mathcal{L}_\tau$ satisfies
\[
\left|x-\frac{p}{q}\right|
\geq q^{-\tau}
\]
for all sufficiently large $q$.
Applying this to its convergents yields
\[
q_{n+1}<q_n^{\tau-1}
\]
for all sufficiently large $n$, which implies
$x\in\mathcal{BR}_\alpha$ by the same geometric summability argument
as above.
Thus
\[
\Delta^\alpha\setminus\mathcal{BR}_\alpha
\subseteq
\bigcap_{\tau=3}^{\infty}\mathcal{L}_\tau.
\]
For every integer $Q\geq2$, the intervals of radius $q^{-\tau}$ about
$p/q$, with $q\geq Q$ and $0\leq p\leq q$, cover
$\mathcal{L}_\tau$.
For any $d>2/\tau$, their $d$th powers of diameters have sum at most
\[
\sum_{q=Q}^{\infty}
(q+1)(2q^{-\tau})^d,
\]
which tends to zero as $Q\to\infty$.
Hence
\[
\dim_H\mathcal{L}_\tau
\leq
\frac{2}{\tau}.
\]
Letting $\tau\to\infty$ proves that
\[
\dim_H
\bigl(
\Delta^\alpha\setminus\mathcal{BR}_\alpha
\bigr)
=0.
\]
\end{proof}

\vspace{1em}

\section{A quantitative analytic iteration}\label{sec:analytic-newton}

For $r>0$, write
\[
 \mathbb T_r:=\{z\in\mathbb C/\mathbb Z:|\operatorname{Im}z|<r\},
 \qquad |g|_r:=\sup_{z\in\mathbb T_r}|g(z)|.
\]
All functions to which $|\cdot|_r$ is applied are periodic and
bounded on the indicated strip. They are called real if they take real values on $\mathbb T$.
We denote the space of bounded holomorphic periodic functions on
$\mathbb T_r$ by $\mathcal A_r$. With the supremum norm it is a
Banach algebra: a uniformly Cauchy sequence converges uniformly,
its limit is holomorphic by local uniform convergence, and
$|gh|_r\le |g|_r|h|_r$. Restriction to a narrower strip has norm
at most one. We will use the following elementary estimates:
\begin{equation}\label{eq:elementary-analytic-estimates}
 |g^{(j)}|_{R-t}\le j!t^{-j}|g|_R,
 \qquad |\langle g\rangle|\le |g|_R,
 \qquad |\widehat g(k)|\le |g|_R e^{-2\pi R|k|}.
\end{equation}
Here $j\ge0$, $0<t<R$, and
$\widehat g(k):=\int_0^1g(x)e^{-2\pi ikx}\,dx$.
The derivative estimate follows from Cauchy's integral formula on
discs of radius less than $t$, followed by a limit. To obtain the
Fourier estimate for $k\ne0$, move the contour to
$\operatorname{Im}z=-\operatorname{sgn}(k)R'$ with $0<R'<R$.
The vertical sides cancel by periodicity; letting $R'\uparrow R$
gives the stated bound. In particular, the Fourier series converges
absolutely on every strictly smaller strip and equals $g$ there.

We write $g_\pm(\theta):=g(\theta\pm\omega)$ and set
\[
 D_+g:=g_+-g,\qquad D_-g:=g-g_-,\qquad
 \Delta_\omega g:=D_-D_+g=g_+-2g+g_-.
\]

\subsection{The conjugacy equation}

For a real periodic forcing $f$ of zero average, consider
\[
 F_f(x,y):=(x+\omega+y+f(x),\ y+f(x)).
\]
A parametrization $K(\theta):=(\theta+u(\theta),v(\theta))$ satisfies
$F_f\circ K=K\circ R_\omega$ precisely when
\begin{equation}\label{eq:scalar-conjugacy}
 v(\theta)=u(\theta)-u(\theta-\omega),\qquad
 u(\theta+\omega)-2u(\theta)+u(\theta-\omega)
       =f(\theta+u(\theta)).
\end{equation}
Indeed, subtracting the second component of the conjugacy identity
from the first gives the formula for $v(\theta+\omega)$, and the
remaining component gives the difference equation. The converse is
verified by substitution. The condition $\langle u\rangle=0$ fixes
the freedom to translate the parameter: if $u$ is a solution, then
$u_c(\theta):=c+u(\theta+c)$ is also a solution and
$\langle u_c\rangle=c+\langle u\rangle$.
Once $\|u'\|_{C^0}<1$,
$\operatorname{Id}+u$ is an orientation-preserving analytic circle
diffeomorphism, so the parametrized curve is a graph.

\subsection{Solving first-order difference equations}

\begin{Lemma}\label{lem:analytic-cohomology}
Suppose that $0<p<1$ and that $D,A>0$ satisfy
\begin{equation}\label{eq:stretched-divisor}
 \|k\omega\|^{-1}\le D\exp(A|k|^p)
 \qquad(k\in\mathbb Z\setminus\{0\}).
\end{equation}
Put $\beta:=p/(1-p)$. There is $B>0$, depending only on $p,D,A$,
such that, whenever $0<t\le1$, $R>t$, and $g\in\mathcal A_R$
with $\langle g\rangle=0$,
each equation $D_\pm w=g$ has a unique zero-average solution on
$\mathbb T_{R-t}$, with
\begin{equation}\label{eq:analytic-cohomology-bound}
 |w|_{R-t}\le\exp(Bt^{-\beta})|g|_R.
\end{equation}
If $g$ is real, then $w$ is real.
\end{Lemma}

\begin{proof}
We consider $D_+w=g$; the equation $D_-w=g$ is treated similarly.
Since $\widehat g(0)=0$, we divide the $k$th Fourier coefficient
of $g$ by $e^{2\pi ik\omega}-1$ for $k\ne0$ and set the
coefficient at zero equal to zero. This gives
\[
w(\theta)=
\sum_{k\ne0}
\frac{\widehat g(k)}
{e^{2\pi ik\omega}-1}
e^{2\pi ik\theta}.
\]
Since
\[
 |e^{2\pi ik\omega}-1|\ge4\|k\omega\|,
 \qquad |\widehat g(k)|\le |g|_R e^{-2\pi R|k|},
\]
we obtain, for either equation,
\[
 |w|_{R-t}\le\frac D2 |g|_R\sum_{k\ge1}
       \exp(Ak^p-2\pi tk).
\]
Maximizing $Ax^p-\pi tx$ over $x\ge0$ gives
\[
 Ax^p-\pi tx\le
 (1-p)A^{1/(1-p)}
 \left(\frac p\pi\right)^{p/(1-p)}t^{-\beta}.
\]
Since $\sum_{k\ge1}e^{-\pi tk}=(e^{\pi t}-1)^{-1}\le(\pi t)^{-1}$,
this gives
\[
 |w|_{R-t}
\le\frac{D}{2\pi t}
 \exp\!\left((1-p)A^{1/(1-p)}(p/\pi)^\beta t^{-\beta}\right)|g|_R.
\]
For $0<t\le1$, elementary maximization gives
$\log(1/t)\le t^{-\beta}/(\beta\exp(1))$.
Consequently one may take
\[
 B:=(1-p)A^{1/(1-p)}(p/\pi)^\beta
       +\max\{0,\log(D/(2\pi))\}
       +\frac{1}{\beta\exp(1)}+1.
\]
This proves \eqref{eq:analytic-cohomology-bound} and also gives
uniform absolute convergence on $\mathbb T_{R-t}$. Substitution
of the Fourier series verifies the equation there. Any two solutions
have a difference whose nonzero Fourier coefficients vanish, since
$\omega$ is irrational; the zero-average normalization removes the
remaining constant. For real $g$, conjugation of the coefficients
shows that $w$ is real.
\end{proof}

\subsection{The analytic existence estimate}

\begin{Proposition}\label{prop:analytic-newton}
Assume \eqref{eq:stretched-divisor} and set $\beta=p/(1-p)$.
There exist constants $C_0>C_1>0$, depending only on $p,D,A$,
with the following property. If $0<r\le1/8$ and $f$ is real
analytic on $\mathbb T_{4r}$, with $\langle f\rangle=0$ and
\begin{equation}\label{eq:newton-smallness}
 |f|_{4r}\le\frac r{16}\exp(-C_0r^{-\beta}),
\end{equation}
then there is a real analytic zero-average function $u$ on
$\mathbb T_r$ satisfying
\begin{equation}\label{eq:newton-solution}
 \Delta_\omega u=f(\theta+u(\theta)).
\end{equation}
Moreover,
\begin{equation}\label{eq:newton-solution-bound}
 |u|_r+|u'|_r
 \le 2\exp(C_1r^{-\beta})|f|_{4r}
 \le\frac r8.
\end{equation}
\end{Proposition}

\begin{proof}
If $f=0$, take $u=0$. Otherwise we give the construction and the
estimates, including the two average conditions needed at each step.

\smallskip\noindent
\textbf{Step 1. Compatibility identity.}
For a real periodic analytic trial function $u$, put
\[
 h:=\operatorname{Id}+u.
\]
Whenever the composition is defined, one has
\begin{equation}\label{eq:newton-exact-compatibility}
 \langle h'[\Delta_\omega u-f(h)]\rangle=0.
\end{equation}
Indeed, on the lifted strip set $V(z):=\int_0^z f(\zeta)\,d\zeta$.
The strip is simply connected, so this defines a holomorphic primitive.
The difference $V(z+1)-V(z)$ is constant, and equals
$\int_0^1 f(x)\,dx=0$. Thus $V$ is periodic, and
$\langle h'f(h)\rangle=\langle(V\circ h)'\rangle=0$.
Also $\langle\Delta_\omega u\rangle=0$. Symmetry of
$\Delta_\omega$ for the bilinear pairing $\langle gh\rangle$,
commutation with differentiation, and integration by parts imply
\[
 \langle u'\Delta_\omega u\rangle
 =\langle u(\Delta_\omega u)'\rangle
 =-\langle u'\Delta_\omega u\rangle=0.
\]
Thus \eqref{eq:newton-exact-compatibility} holds for every trial
function, not only for a solution.

\smallskip\noindent
\textbf{Step 2. Construction and estimate of a correction.}
Let $r<R\le2r$, $0<t\le r$, and $R-4t\ge r$.
Suppose that $u$ is real, $\langle u\rangle=0$,
$|u|_R\le r/2$, and $|u'|_R\le1/4$.
Then $h(\mathbb T_R)\subset\mathbb T_{5r/2}$, so
$\Delta_\omega u-f(h)$ belongs to $\mathcal A_R$.
The linearization of the equation is
\[
 Lv:=\Delta_\omega v-f'(h)v.
\]
Differentiation of the residual gives $Lh'=\partial_\theta\bigl[\Delta_\omega u-f(h)\bigr]$.
Then we have
\[
 D_-(h'h'_+D_+w)
 =h'\bigl[h'_+(w_+-w)+h'_-(w_--w)\bigr]
 =h'\bigl[\Delta_\omega(h'w)-w\Delta_\omega h'\bigr].
\]
Substituting $Lh'=\partial_\theta\bigl[\Delta_\omega u-f(h)\bigr]$ gives the exact identity
\begin{equation}\label{eq:newton-factorization}
 L(h'w)=\frac{1}{h'}D_-\bigl(h'h'_+D_+w\bigr)+w\partial_\theta\bigl[\Delta_\omega u-f(h)\bigr].
\end{equation}
We solve two first-order equations. First let $z$ be the
zero-average solution of
\[
 D_-z=-h'[\Delta_\omega u-f(h)].
\]
It exists by \eqref{eq:newton-exact-compatibility} and
Lemma~\ref{lem:analytic-cohomology}. On the real circle,
$3/4\le h'\le5/4$, and consequently
\[
 \langle \frac{1}{h'h'_+}\rangle\ge16/25.
\]
On the complex strip one has $|\frac{1}{h'h'_+}|_R\le16/9$.
Define
\[
 c:=-\frac{\langle \frac{z}{h'h'_+}\rangle}{\langle \frac{1}{h'h'_+}\rangle},
 \qquad D_+w_0=\frac{1}{h'h'_+}(z+c),\qquad \langle w_0\rangle=0.
\]
The definition of $c$ gives the compatibility condition for the
second equation. Finally, set
\[
 d:=-\langle h'w_0\rangle,\qquad w:=w_0+d,\qquad v:=h'w.
\]
Since $\langle h'\rangle=1$, this choice gives $\langle v\rangle
=\langle h'w_0\rangle+d\langle h'\rangle
=0$.
The solutions $z,w_0$ and the constants $c,d$ are real.
Equation~\eqref{eq:newton-factorization} yields
\begin{equation}\label{eq:newton-linear-defect}
 Lv=-[\Delta_\omega u-f(h)]+w\partial_\theta\bigl[\Delta_\omega u-f(h)\bigr].
\end{equation}

The first inversion and the bounds for $h'$ give
\[
 |z|_{R-t}\le\frac54e^{Bt^{-\beta}}|\Delta_\omega u-f(h)|_R.
\]
Since $\langle \frac{1}{h'h'_+}\rangle^{-1}\le25/16$ and
$|\frac{1}{h'h'_+}|_R\le16/9$, the constant $c$ and the right-hand side of
the second equation satisfy
\[
 |c|\le\frac{25}{9}|z|_{R-t},\qquad
 |\frac{(z+c)}{h'h'_+}|_{R-t}\le\frac{544}{81}|z|_{R-t}.
\]
Applying Lemma~\ref{lem:analytic-cohomology} a second time yields
\[
 |w_0|_{R-2t}\le\frac{680}{81}e^{2Bt^{-\beta}}|\Delta_\omega u-f(h)|_R.
\]
Furthermore, $|d|\le(5/4)|w_0|_{R-2t}$ and $v=h'w$, so
\begin{equation}\label{eq:newton-correction-estimate}
 \begin{split}
 |w|_{R-2t}&\le\frac{170}{9}e^{2Bt^{-\beta}}|\Delta_\omega u-f(h)|_R\le20 e^{2Bt^{-\beta}}|\Delta_\omega u-f(h)|_R,\\
 |v|_{R-2t}&\le\frac{425}{18}e^{2Bt^{-\beta}}|\Delta_\omega u-f(h)|_R\le25e^{2Bt^{-\beta}}|\Delta_\omega u-f(h)|_R,\\
 |v'|_{R-3t}&\le25t^{-1}e^{2Bt^{-\beta}}|\Delta_\omega u-f(h)|_R.
 \end{split}
\end{equation}
The last inequality is Cauchy's estimate applied directly to $v$.
In particular, it does not require a bound for $u''$.

If $|v|_{R-2t}\le r/2$, then, for $z\in\mathbb T_{R-2t}$
and $0\le s\le1$,
\[
 |\operatorname{Im}(h(z)+sv(z))|
 \le |\operatorname{Im}z|+|u(z)|+|v(z)|<3r.
\]
Thus the entire line segment used in Taylor's formula lies in
$\mathbb T_{3r}$. Identity
\eqref{eq:newton-linear-defect} gives
\begin{equation}\label{eq:newton-quadratic-residual}
 \Delta_\omega(u+v)-f(h+v)=w\partial_\theta[\Delta_\omega u-f(h)]-v^2\int_0^1(1-s)f''(h+sv)\,ds.
\end{equation}
Cauchy's estimates imply
$|\partial_\theta[\Delta_\omega u-f(h)]|_{R-t}\le t^{-1}|\Delta_\omega u-f(h)|_R$ and
$|f''|_{3r}\le2r^{-2}|f|_{4r}$.
Since \eqref{eq:newton-smallness} implies $|f|_{4r}\le1$,
and $t\le r\le1$, we obtain
\[
 \begin{split}
 |\Delta_\omega(u+v)-f(h+v)|_{R-4t}
 &\le\bigl(20t^{-1}e^{2Bt^{-\beta}}
       +625r^{-2}|f|_{4r}e^{4Bt^{-\beta}}\bigr)|\Delta_\omega u-f(h)|_R^2\\
 &\le700t^{-2} e^{4Bt^{-\beta}}|\Delta_\omega u-f(h)|_R^2.
 \end{split}
\]
For $0<t\le1$, the inequality
$\log(1/t)\le t^{-\beta}/(\beta\exp(1))$ shows that one may take
\[
 B_1:=4B+\log700+\frac{2}{\beta\exp(1)}.
\]
Indeed, this choice gives both
\begin{align}
 |v|_{R-2t}+|v'|_{R-3t}
 &\le e^{B_1t^{-\beta}}|\Delta_\omega u-f(h)|_R,
 \label{eq:newton-one-step-correction}\\
 |\Delta_\omega(u+v)-f(h+v)|_{R-4t}
 &\le e^{B_1t^{-\beta}}|\Delta_\omega u-f(h)|_R^2.
 \label{eq:newton-one-step-error}
\end{align}
The first estimate is established before, and independently of,
the condition $|v|_{R-2t}\le r/2$ needed for the second estimate.

\smallskip\noindent
\textbf{Step 3. Convergence.}

Choose
\[
\tau:=2^{-1/(2\beta)},\qquad
R_j:=r+r\tau^j,\qquad
t_j:=\frac{R_j-R_{j+1}}{4},
\]
and set
\[
C_1:=B_1\left(\frac{4}{1-\tau}\right)^\beta,
\qquad
C_0:=\frac{2C_1}{2-\sqrt{2}}.
\]
Then
\[
R_{j+1}=R_j-4t_j,
\qquad
B_1t_j^{-\beta}=C_1r^{-\beta}2^{j/2}.
\]
The smallness assumption gives
\[
e^{(C_0/2)r^{-\beta}}|f|_{2r}
\le
\frac{r}{16}e^{-(C_0/2)r^{-\beta}}
<\frac12,
\qquad
2e^{C_1r^{-\beta}}|f|_{2r}\le\frac r8.
\]

Start with $u_0=0$, and write
$h_j:=\operatorname{Id}+u_j$ whenever $u_j$ has been constructed.
We construct $u_{j+1}=u_j+v_j$ inductively and prove that
\begin{equation}\label{eq:direct-iteration-residual}
|\Delta_\omega u_j-f(h_j)|_{R_j}
\le
|f|_{2r}^{\,2^j}
\exp\!\left[
\frac{C_0}{2}r^{-\beta}
\left(2^j-2^{j/2}\right)
\right]
\end{equation}
and
\begin{equation}\label{eq:direct-iteration-correction}
|v_j|_{R_j-2t_j}+|v_j'|_{R_j-3t_j}
\le
2^{-j}e^{C_1r^{-\beta}}|f|_{2r}.
\end{equation}

The residual estimate holds for $j=0$.
Suppose that $u_j$ has been constructed,
\eqref{eq:direct-iteration-residual} holds at step $j$,
and \eqref{eq:direct-iteration-correction} holds for all $i<j$.
Restricting the previous corrections to $\mathbb T_{R_j}$ gives
\[
|u_j|_{R_j}+|u_j'|_{R_j}
\le
e^{C_1r^{-\beta}}|f|_{2r}
\sum_{i<j}2^{-i}
\le\frac r8.
\]
Thus the two first-order inversions in Step~2 construct a real
zero-average correction $v_j$.

Since $C_0/2>C_1$, the correction estimate and
\eqref{eq:direct-iteration-residual} imply
\[
\begin{aligned}
|v_j|_{R_j-2t_j}+|v_j'|_{R_j-3t_j}
&\le
e^{C_1r^{-\beta}2^{j/2}}
|\Delta_\omega u_j-f(h_j)|_{R_j}\\
&\le
e^{C_1r^{-\beta}}|f|_{2r}
\left(
e^{(C_0/2)r^{-\beta}}|f|_{2r}
\right)^{2^j-1}\\
&\le
2^{-j}e^{C_1r^{-\beta}}|f|_{2r},
\end{aligned}
\]
where the last inequality uses $2^j-1\ge j$.
This proves \eqref{eq:direct-iteration-correction} at step $j$.
In particular,
\[
|v_j|_{R_j-2t_j}\le\frac r{16}<\frac r2.
\]
Define $u_{j+1}:=u_j+v_j$.
The Taylor segment therefore lies in the domain of $f$, so the
residual estimate from Step~2 applies.

Using
\[
C_1=\frac{C_0}{2}(2-\sqrt{2}),
\]
we obtain
\[
\begin{aligned}
|\Delta_\omega u_{j+1}-f(h_{j+1})|_{R_{j+1}}
&\le
e^{C_1r^{-\beta}2^{j/2}}
|\Delta_\omega u_j-f(h_j)|_{R_j}^{\,2}\\
&\le
|f|_{2r}^{\,2^{j+1}}
\exp\!\left[
\frac{C_0}{2}r^{-\beta}
\left(2^{j+1}-2^{(j+1)/2}\right)
\right].
\end{aligned}
\]
This is \eqref{eq:direct-iteration-residual} at step $j+1$,
and the induction is complete.

Summing \eqref{eq:direct-iteration-correction} yields
\[
\sum_{j\ge0}
\left(
|v_j|_{R_j-2t_j}+|v_j'|_{R_j-3t_j}
\right)
\le
2e^{C_1r^{-\beta}}|f|_{2r}
\le\frac r8.
\]
Consequently, $\sum_j v_j$ and $\sum_j v_j'$ converge uniformly
on $\mathbb T_r$.
The first series defines a bounded holomorphic function $u$ there,
and local Cauchy integral formulas identify the second limit with $u'$.
Real-valuedness and zero average pass to the limit, and
\[
|u|_r+|u'|_r
\le
2e^{C_1r^{-\beta}}|f|_{2r}
\le
2e^{C_1r^{-\beta}}|f|_{4r}
\le\frac r8.
\]

Moreover, \eqref{eq:direct-iteration-residual} gives
\[
\begin{aligned}
|\Delta_\omega u_j-f(h_j)|_{R_j}
&\le
\left(
e^{(C_0/2)r^{-\beta}}|f|_{2r}
\right)^{2^j}
\exp\!\left[
-\frac{C_0}{2}r^{-\beta}2^{j/2}
\right]\\
&\le 2^{-2^j}\longrightarrow0.
\end{aligned}
\]
The maps $h_j$ and $h:=\operatorname{Id}+u$ send
$\mathbb T_r$ into $\mathbb T_{9r/8}$.
Cauchy's estimate, followed by integration along a line segment, gives
\[
|f(h_j)-f(h)|_r
\le
r^{-1}|f|_{4r}|u_j-u|_r
\longrightarrow0.
\]
Also,
\[
|\Delta_\omega(u_j-u)|_r
\le4|u_j-u|_r
\longrightarrow0.
\]
Passing to the limit proves
\[
\Delta_\omega u=f(\operatorname{Id}+u)
\]
on $\mathbb T_r$.
This completes the proof of the proposition.
\end{proof}

\vspace{1em}

\section{Proof of the main theorem}\label{sec:proofth1}

\subsection{Arithmetic bounds for the difference equation}
\label{sec:newton-arithmetic}

The arithmetic condition needed below is the pointwise condition
\eqref{he}, which is equivalent to $\omega\in\Delta^\alpha$ by
Proposition~\ref{bca}. Fix $H\ge1$ for which that condition holds.
For every $j\in\mathbb Z\setminus\{0\}$, choose the largest index $n$
such that $q_n\le |j|$. Proposition~\ref{prop:continued-fraction-estimates} and \eqref{he} give
\begin{equation}\label{eq:newton-small-divisor}
 \|j\omega\|_{\mathbb R/\mathbb Z}
 \ge |q_n\omega-p_n|
 >\frac1{q_n+q_{n+1}}
 \ge\frac12\exp\bigl(-H|j|^{1/\alpha}\bigr).
\end{equation}
Here $q_n\le q_{n+1}\le e^{Hq_n^{1/\alpha}}$ was used.
Choosing the largest index ensures $|j|<q_{n+1}$ and also handles
the possible repetition $q_0=q_1=1$. Thus \eqref{eq:stretched-divisor} holds with $p=1/\alpha$,
$D=2$ and $A=H$, and the exponent in
Proposition~\ref{prop:analytic-newton} is $\beta=1/(\alpha-1)$.
Its constants depend on the fixed $\alpha$ and $\omega$, and are
independent of $m$, $N_m$, and the choice of the sequence $(N_m)$.
The estimate \eqref{eq:newton-small-divisor} holds for every
nonzero Fourier index, including every multiple of $N_m$.
Throughout the iteration the difference operators are
$D_\pm$ at the fixed rotation number $\omega$; their inverses
are estimated on the original angle variable $\theta$. Thus no
additional arithmetic condition on $N_m\omega$ is needed.

\subsection{Analytic size of the rescaled perturbation}
\label{sec:newton-rescaled-force}

Suppose first that $0<\gamma<1$ and
$\varphi\in G^{\gamma,R_\varphi}(\mathbb T)$ for some $R_\varphi>0$.
There are constants $C_\varphi,a_\varphi>0$ such that
\begin{equation}\label{eq:newton-ultra-fourier}
 |\widehat\varphi(j)|
 \le C_\varphi\exp\bigl(-a_\varphi|j|^{1/\gamma}\bigr),
 \qquad j\in\mathbb Z.
\end{equation}
The defining norm gives
\[
 |\widehat\varphi(j)|
 \le \|\varphi\|_{\gamma,R_\varphi}
       \frac{R_\varphi^{-\gamma k}(k!)^\gamma}{(2\pi |j|)^k}
 \qquad(j\ne0,\ k\ge0).
\]
Take $k=\lfloor R_\varphi(2\pi|j|)^{1/\gamma}/2\rfloor$ and use
$k!\le k^k$. For large $|j|$ the resulting bound is at most
$\|\varphi\|_{\gamma,R_\varphi}2^{-\gamma k}$, which proves
\eqref{eq:newton-ultra-fourier}; a larger constant covers the finitely
many remaining indices. More explicitly, once
$R_\varphi(2\pi|j|)^{1/\gamma}\ge4$, this choice has
$k\ge R_\varphi(2\pi|j|)^{1/\gamma}/4$, so one may take
$a_\varphi=(\gamma\log2)R_\varphi(2\pi)^{1/\gamma}/4$.
Since $1/\gamma>1$, the series
\[
 \sum_{j\in\mathbb Z}\widehat\varphi(j)e^{2\pi ijz}
\]
and all its termwise derivatives converge uniformly on every finite
closed strip. For each fixed $\ell\ge0$ and $S>0$,
\[
 |j|^\ell e^{2\pi S|j|}
 \le e^{(a_\varphi/2)|j|^{1/\gamma}}
\]
for all sufficiently large $|j|$, so \eqref{eq:newton-ultra-fourier}
gives a summable bound for each differentiated series. This defines an entire periodic function. On the real circle
it equals $\varphi$, because the Fourier series of the smooth
function $\varphi$ converges uniformly to it. We use this extension
when taking analytic strip norms.

Let $N\ge1$ be an integer and $0<s\le1/8$. Splitting the Fourier decay
in \eqref{eq:newton-ultra-fourier} into two equal parts yields
\begin{equation}\label{eq:newton-force-width}
 \begin{aligned}
 |\varphi(N\cdot)|_{4s}
 &\le C_\varphi
   \sum_{j\in\mathbb Z}
    e^{-a_\varphi|j|^{1/\gamma}+8\pi Ns|j|}\\
 &\le F_0\exp\left(F_1(Ns)^{1/(1-\gamma)}\right).
 \end{aligned}
\end{equation}
For $c\ge0$ and $d>0$, differentiation gives
\[
 \sup_{x\ge0}\{cx-dx^{1/\gamma}\}
 =(1-\gamma)(\gamma/d)^{\gamma/(1-\gamma)}c^{1/(1-\gamma)}.
\]
Taking $c=8\pi Ns$ and $d=a_\varphi/2$ proves the second
inequality in \eqref{eq:newton-force-width}, with
\[
 \begin{split}
 F_0&:=C_\varphi\sum_{j\in\mathbb Z}
                  e^{-(a_\varphi/2)|j|^{1/\gamma}},\\
 F_1&:=(1-\gamma)(2\gamma/a_\varphi)^{\gamma/(1-\gamma)}
                           (8\pi)^{1/(1-\gamma)}.
 \end{split}
\]
These constants depend on the fixed $\varphi,\gamma,R_\varphi$,
but not on $N$ or $s$.
If $\gamma=0$, let $J_\varphi$ be the Fourier degree of $\varphi$ and write
$\varphi=\sum_{|j|\le J_\varphi}\widehat\varphi(j)e^{2\pi ij\theta}$.
Then
\[
 |\varphi(N\cdot)|_{4s}
 \le\left(\sum_{|j|\le J_\varphi}|\widehat\varphi(j)|\right)e^{8\pi J_\varphi Ns}.
\]
Thus \eqref{eq:newton-force-width} holds also in this case, with
$F_0:=\sum_{|j|\le J_\varphi}|\widehat\varphi(j)|$, $F_1:=8\pi J_\varphi$, and
$1/(1-\gamma)=1$. Here $J_\varphi\ge1$, since $\varphi$ is nonzero and
has zero average. No zero-index derivative norm is used.

\subsection{Choice of strip width and perturbation size}
\label{sec:newton-main-application}

Set
\begin{equation}\label{eq:newton-optimal-exponents}
 \sigma:=\frac{\alpha-1}{\alpha-\gamma},\qquad
 \kappa:=\frac1{\alpha-\gamma},\qquad
 s_m:=\tfrac18N_m^{-\sigma}.
\end{equation}
Since $0\le\gamma<1$, we have $0<\sigma<1$. The two
exponents in the analytic estimates have the same order:
\[
 s_m^{-1/(\alpha-1)}=8^{1/(\alpha-1)}N_m^\kappa,
 \qquad
 (N_ms_m)^{1/(1-\gamma)}=8^{-1/(1-\gamma)}N_m^\kappa.
\]
The hypothesis $\gamma<\varepsilon/\alpha$ in
Theorem~\ref{Th1} implies $\gamma<\varepsilon$, so
\[
 \frac1{\alpha-\gamma}<\frac1{\alpha-\varepsilon}.
\]
Use the exponent and amplitudes in Theorem~\ref{Th1}:
\begin{equation}\label{eq:newton-delta}
 \frac1{\alpha-\gamma}<\tau_*<\frac1{\alpha-\varepsilon},
 \qquad \delta_m=\exp(-N_m^{\tau_*}).
\end{equation}
The choice is independent of $L$. For the function $\phi_m$ defined in the introduction,
\eqref{eq:newton-force-width} gives
\begin{equation}\label{eq:newton-initial-residual}
 \delta_m|\phi_m|_{4s_m}
 \le F_0\exp\left(-N_m^{\tau_*}+F_1N_m^\kappa\right).
\end{equation}
The smallness threshold in Proposition~\ref{prop:analytic-newton}
is $(s_m/16)\exp(-C_0s_m^{-\beta})$, where
$\beta=1/(\alpha-1)$ and $C_0$ is fixed. In particular,
\begin{equation}\label{eq:threshold-ratio}
 \begin{split}
 \log\frac{\delta_m|\phi_m|_{4s_m}}
                   {(s_m/16)\exp(-C_0s_m^{-\beta})}
 &\le -N_m^{\tau_*}
       +(F_1+C_0 8^\beta)N_m^\kappa\\
 &\qquad +\sigma\log N_m+\log(128F_0)
 \longrightarrow-\infty.
 \end{split}
\end{equation}
This follows from $\tau_*>\kappa$. To make the uniformity in the
oscillation frequency explicit, use
$\log N\le N^\kappa/(\kappa\exp(1))$ for every real $N\ge1$ and set
\[
 C_*:=F_1+C_0 8^\beta
       +\frac{\sigma}{\kappa\exp(1)}
       +\max\{0,\log(128F_0)\}.
\]
The right-hand side of \eqref{eq:threshold-ratio}, with $N$ in
place of $N_m$, is at most $-N^{\tau_*}+C_*N^\kappa$.
Consequently we may choose
\begin{equation}\label{eq:uniform-frequency-threshold}
 N_*:=\left\lceil
       \max\{1,(2C_*)^{1/(\tau_*-\kappa)}\}
       \right\rceil.
\end{equation}
For every integer $N\ge N_*$, this bound is at most
$-\tfrac12N^{\tau_*}<0$. Thus the smallness condition holds
whenever $N_m\ge N_*$. The threshold depends only on the fixed
$\alpha,\omega,\gamma,\varphi,R_\varphi,\tau_*$ and is independent
of $L$ and of the sequence $(N_m)$. Here $R_\varphi$ is chosen
once for the fixed function $\varphi$ and is omitted when
$\gamma=0$. Since $N_m\to\infty$, there is an index $m_0$ such
that $N_m\ge N_*$ for every $m\ge m_0$; this index may depend
on the chosen sequence.

The forcing $\delta_m\phi_m$ is real and has zero average by
the integer rescaling identity in the introduction. Its derivatives
are controlled on the same angle strips: Cauchy's estimate gives
\[
 \delta_m|\phi_m^{(\ell)}|_{3s_m}
 \le \ell!s_m^{-\ell}\delta_m|\phi_m|_{4s_m},
 \qquad \ell\ge0.
\]
In particular, the derivative factors arising from
$\varphi(N_m(\theta+u(\theta)))$ are already controlled by
the strip norm used in Proposition~\ref{prop:analytic-newton}.
That proposition, applied at the fixed rotation number $\omega$,
provides a real-analytic zero-average solution $u_m$ of
\[
 u_m(\theta+\omega)-2u_m(\theta)+u_m(\theta-\omega)
   =\delta_m\varphi\bigl(N_m(\theta+u_m(\theta))\bigr).
\]
The solution is analytic on a strip of width at least $s_m$.
Its derivative satisfies
\begin{equation}\label{eq:newton-graph-derivative}
 \|u_m'\|_{C^0}
 \le F_2\exp\left(-N_m^{\tau_*}+F_3N_m^\kappa\right)
 \longrightarrow0.
\end{equation}
For instance, \eqref{eq:newton-solution-bound} and
\eqref{eq:newton-initial-residual} permit
$F_2:=2F_0$ and $F_3:=F_1+C_1 8^\beta$.
The constants $F_0,F_1$ have the dependencies specified above;
$F_2,F_3$ may also depend on $\alpha$ and $\omega$.
None of them depends on $m$, on the sequence $(N_m)$, or on $L$.

The proposition also gives the complex-strip estimate
\begin{equation}\label{eq:graph-strip-smallness}
 |u_m|_{s_m}+|u_m'|_{s_m}\le s_m/8.
\end{equation}
The lift $h_m(\theta):=\theta+u_m(\theta)$ has
positive derivative on the real line and satisfies
$h_m(\theta+1)=h_m(\theta)+1$. It is bijective there and induces
an orientation-preserving circle diffeomorphism.

We construct its analytic inverse explicitly. Set $r=s_m$ and use
lifts to the complex plane. For each $x$ with
$|\operatorname{Im}x|<r/2$, the map
\[
 z\longmapsto x-u_m(z)
\]
sends the closed disc $\overline D(x,r/4)$ into
$\overline D(x,r/8)$ and has Lipschitz constant at most $r/8<1$.
The disc lies in $\mathbb T_{3r/4}$, so these assertions follow
from \eqref{eq:graph-strip-smallness} and integration of $u_m'$
along a line segment. Starting at $z_0(x)=x$, the iteration
$z_{\ell+1}(x)=x-u_m(z_\ell(x))$ converges uniformly for
$|\operatorname{Im}x|<r/2$. Indeed,
\[
 |z_{\ell+1}(x)-z_\ell(x)|\le(r/8)^{\ell+1}.
\]
Every $z_\ell$ is holomorphic, so its limit $g_m$ is holomorphic.
It satisfies $h_m(g_m(x))=x$ and $|g_m(x)-x|\le r/8$;
thus $g_m(\mathbb T_{r/2})\subset\mathbb T_{5r/8}$.
Uniqueness of the fixed point gives
$g_m(x+1)=g_m(x)+1$ and real-valuedness on the real line.
Consequently its real restriction is $h_m^{-1}$.

Define
\[
 v_m(\theta):=u_m(\theta)-u_m(\theta-\omega),\qquad
 \Psi_m:=v_m\circ g_m.
\]
The function $\Psi_m$ is real, one-periodic, and holomorphic on
$\mathbb T_{s_m/2}$. The scalar equation implies, component by
component,
\[
 \begin{split}
 v_m(\theta)+\delta_m\phi_m(h_m(\theta))
     &=u_m(\theta+\omega)-u_m(\theta)=v_m(\theta+\omega),\\
 h_m(\theta)+\omega+v_m(\theta)+\delta_m\phi_m(h_m(\theta))
     &=h_m(\theta+\omega).
 \end{split}
\]
Therefore
$f_m(h_m(\theta),v_m(\theta))=(h_m(\theta+\omega),v_m(\theta+\omega))$.
As $h_m$ parametrizes the entire circle, the graph of $\Psi_m$
has Property $A_\omega$.

\subsection{The two Gevrey norms}
\label{sec:newton-gevrey-asymptotics}

For each positive integer $N_m$, the map $x\mapsto N_mx$ is
onto $\T$. The chain rule consequently gives the exact identity
\begin{equation}\label{eq:integer-rescaling-derivatives}
 \|\phi_m^{(k)}\|_{C^0(\T)}
 =N_m^k\|\varphi^{(k)}\|_{C^0(\T)},
 \qquad k\ge0.
\end{equation}
This identity applies to both the upper and lower bounds below.
Fix $L>0$. When $0<\gamma<1$, choose $R_\varphi>0$ such that
$\varphi\in G^{\gamma,R_\varphi}$. Each term in the defining norm gives
\[
 \|\varphi^{(k)}\|_{C^0}
 \le \|\varphi\|_{\gamma,R_\varphi}R_\varphi^{-\gamma k}(k!)^\gamma.
\]
For any $a\ge1$ and $X\ge0$, one has
\[
 \sum_{k\ge0}\frac{X^k}{(k!)^a}
 =\sum_{k\ge0}\left(\frac{X^{k/a}}{k!}\right)^a
 \le\left(\sum_{k\ge0}\frac{X^{k/a}}{k!}\right)^a
 =\exp(aX^{1/a}).
\]
For the inequality, divide a nonnegative sequence by its sum; each
resulting term lies in $[0,1]$ and its $a$th power is no larger
than itself. Applying this with $a=\alpha-\gamma>1$ gives
\begin{equation}\label{eq:newton-weak-gevrey}
 \begin{aligned}
 \|f_m-T_\omega\|_{\alpha,L}
 &=\delta_m\sum_{k\ge0}
       \frac{L^{\alpha k}N_m^k}{(k!)^\alpha}
           \|\varphi^{(k)}\|_{C^0}\\
 &\le\|\varphi\|_{\gamma,R_\varphi}
   \exp\left[-N_m^{\tau_*}
      +(\alpha-\gamma)
        (L^\alpha R_\varphi^{-\gamma}N_m)^{1/(\alpha-\gamma)}\right]
   \longrightarrow0.
 \end{aligned}
\end{equation}
For $\gamma=0$, let $J_\varphi$ be the Fourier degree fixed in
Subsection~\ref{sec:newton-rescaled-force}. Use
$\|\varphi^{(k)}\|_{C^0}\le F_0(2\pi J_\varphi)^k$ instead.
The same argument with the exponent $\alpha$ gives an upper bound
\[
 F_0\exp\!\left(-N_m^{\tau_*}
       +\alpha(2\pi J_\varphi L^\alpha N_m)^{1/\alpha}\right),
\]
which also tends to zero, since $\tau_*>1/\alpha$.

For the opposite limit, the zero average and nonvanishing of
$\varphi$ imply that
$\widehat\varphi(j_*)\ne0$ for some nonzero integer $j_*$.
Taking this Fourier coefficient in $\varphi^{(k)}$ gives
\[
 \|\varphi^{(k)}\|_{C^0}
 \ge |\widehat\varphi(j_*)|(2\pi|j_*|)^k.
\]
Let $\alpha_-:=\alpha-\varepsilon>1$ and
$B_m:=2\pi|j_*|L^{\alpha_-} N_m$. Then
\begin{equation}\label{eq:newton-strong-gevrey}
 \|f_m-T_\omega\|_{\alpha_-,L}
 \ge |\widehat\varphi(j_*)|e^{-N_m^{\tau_*}}
       \sum_{k\ge0}\frac{B_m^k}{(k!)^{\alpha_-}}.
\end{equation}
For $k\ge1$, monotonicity of $\log x$ gives the elementary
factorial estimate
\[
 \int_1^k\log x\,dx\le\log(k!)
 \le\int_1^k\log x\,dx+\log k
 =k\log k-k+1+\log k.
\]
Put $X_m:=B_m^{1/\alpha_-}$ and $k_m:=\lfloor X_m\rfloor$.
For all sufficiently large $m$, $X_m\ge2$, and the preceding
upper bound for $\log(k_m!)$ implies
\[
 \begin{split}
 \log\frac{B_m^{k_m}}{(k_m!)^{\alpha_-}}
 &\ge\alpha_-\bigl[k_m\log(X_m/k_m)+k_m-1-\log k_m\bigr]\\
 &\ge\alpha_-X_m-\alpha_-(2+\log X_m).
 \end{split}
\]
The last step uses $k_m\le X_m<k_m+1$ and
$\log(X_m/k_m)\ge0$. Keeping just this summand in
\eqref{eq:newton-strong-gevrey} gives
\[
 \begin{split}
 \log\|f_m-T_\omega\|_{\alpha-\varepsilon,L}
 &\ge -N_m^{\tau_*}
       +\alpha_-(2\pi|j_*|L^{\alpha_-})^{1/\alpha_-}N_m^{1/\alpha_-}\\
 &\qquad -\alpha_-(2+\log X_m)+\log|\widehat\varphi(j_*)|
 \longrightarrow+\infty.
 \end{split}
\]
Indeed, $\tau_*<1/\alpha_-$, the coefficient of
$N_m^{1/\alpha_-}$ is positive for every fixed $L>0$, and
$\log X_m=O(\log N_m)$.

Both norms are finite for every fixed $m$. For $\gamma>0$, the
norm of order $s\in\{\alpha,\alpha_-\}$ is bounded by
\[
 \delta_m\|\varphi\|_{\gamma,R_\varphi}
 \sum_{k\ge0}\frac{(L^sR_\varphi^{-\gamma}N_m)^k}{(k!)^{s-\gamma}}.
\]
The ratio of consecutive summands is
$L^sR_\varphi^{-\gamma}N_m/(k+1)^{s-\gamma}$ and tends to zero,
since $s>1>\gamma$. For $\gamma=0$ the finite Fourier estimate
gives the same conclusion with exponent $s>0$ in the factorial.
Thus the limits compare finite norms, and all constants depending
on $L$ enter only after the amplitudes $\delta_m$ and the invariant
graphs have been constructed. The same sequence works for every
fixed $L>0$. Every estimate above depends on the sequence only
through the value of $N_m$ at that step. Thus the two norm limits
hold for every sequence of positive integers tending to infinity,
without a monotonicity or growth-rate assumption. This proves
Theorem~\ref{Th1}.

\vspace{2em}

\noindent\textbf{Data Availability Statement.}
The authors state that this manuscript has no associated data and there is no conflict of interest.

%	 \vspace{2em}
%		
 \section{AI declaration}
The authors used OpenAI’s ChatGPT as an auxiliary tool for proofreading, improving the presentation of the proof of Proposition 3.2. 
All mathematical arguments were  verified by the authors, who take full
responsibility for the content.

\end{document}